\documentclass[12pt,psamsfonts]{amsart} \usepackage{amsmath,amsthm,amsfonts,amssymb}

\usepackage{eucal}

\usepackage{color}
\usepackage{pb-diagram,graphicx, latexsym,enumerate,epsf, epsfig}
\usepackage{lineno}
\usepackage{pstricks}
\usepackage[all,knot,poly]{xy}

 \numberwithin{equation}{section} \newtheorem*{thm}{Theorem}
\newtheorem{theorem}[equation]{Theorem}

\newtheorem{lemma}[equation]{Lemma}

\newtheorem{prop}[equation]{Proposition}
 \theoremstyle{definition}

\newtheorem{rmk}[equation]{Remark}
\newtheorem{remark}[equation]{Remark}

 \DeclareMathOperator{\Ker}{Ker}

\begin{document} \title[Braid Group Representations arising from the Yang Baxter Equation] {Braid Group representations arising from the Yang Baxter Equation}


\author{Jennifer M. Franko} \email{frankoj2@scranton.edu} \address{Department of Mathematics\\
    The University of Scranton \\
    Scranton, PA 18508\\
    U.S.A.}

\thanks{The author was partially supported by NSF FRG grant DMS-034772. and would like
to thank Z. Wang and E. Rowell for helpful correspondence and
support.  This work is based on part of the author's Ph.D. thesis.}

\begin{abstract}

This paper aims to determine the images of the braid group under
representations afforded by the Yang Baxter equation when the
solution is a nontrivial $4 \times 4$ matrix. Making the assumption
that all the eigenvalues of the Yang Baxter solution are roots of
unity, leads to the conclusion that all the images are finite.

 Using results of Turaev,
we have also identified cases in which one would get a link
invariant. Finally, by observing the group algebra generated by the
image of the braid group sometimes factor through known algebras, in
certain instances we can identify the invariant as particular
specializations of a known invariant.
\end{abstract} \maketitle

\section{Introduction} \label{s:intro}

      Any invertible matrix which satisfies the Yang Baxter equation
 can be used to obtain representations of the braid group.  This paper focuses on
cases where the Yang Baxter solution is a $4 \times 4$ unitary
matrix. One of the goals is to identify the image of the braid group
under these representations.  Our main technique is to study the
pure braid group, $P_n$, a subgroup of the braid group, and apply
the results of this analysis to the entire group.  We show that
restricting the eigenvalues of the solution results in the image of
the braid group being finite.

\subsection{}
The $n$-strand braid group, $B_n$, (for $n \geq 2$) is generated by
$\sigma_1,~ \sigma_2,~ ...,~\sigma_{n-1}$ with the following
relations:
\begin{description}
                                                  \item[(B1)] $\sigma_i
                                                  \sigma_j =
                                                  \sigma_j
                                                  \sigma_i~,~|i-j| \geq 2$
                                                  \item[(B2)] $\sigma_i
                                                  \sigma_{i+1}\sigma_i=\sigma_{i+1}
                                                  \sigma_i
                                                  \sigma_{i+1}~,~1
                                                  \leq i \leq n-2$
                                                \end{description}
\indent Representations of the braid group often give rise to knot
and link invariants and Turaev \cite{Turaev} defined a criteria,
called an \emph{enhancement}, which if satisfied, would produce a
Markov trace and hence lead to link invariants.  One approach to
obtaining representations of the braid group is through
consideration of $R$ matrices of quantum groups. However, all that
is required is a solution to the Yang Baxter equation, which can be
systematically produced using the theory of quantum groups.  But
simple solutions of the Yang Baxter equation can also found by
directly solving the equation with the help of computers. Starting
with a finite-dimensional vector space $V$, let $R$ be a linear map
on the tensor product of $V$ with itself, then $R$ is said to
satisfy the Yang Baxter equation if:
$$(R \otimes I)(I \otimes R)(R \otimes I)=(I \otimes R)(R \otimes
I)(I \otimes R).$$ If $R$ satisfies the Yang Baxter equation and if
it is an invertible linear map, it can be used to obtain a
representation of the braid group as follows: $\pi
(\sigma_i)=I^{\otimes (i-1)} \otimes R \otimes I^{\otimes (n-i-1)}$
where $I$ is the identity on $V$. Pictorially, if
$$R=
\begin{picture}(20,30) (20,-10) \thicklines \put( 20,20){\line( 1,-2){
7}} \put(30,0){\line( 1,-2){7}} \put( 37,19){\line(-1,-2){17}}
\end{picture} $$
 this corresponds to the (B2) relation in the braid
group (and the third Reidemiester move).\\








All solutions of the form $R : V \otimes V \rightarrow V \otimes V$
for $V$ of $\text{dimension}~=~2$ to the Yang Baxter equation has
been listed in \cite{H}. Dye found all unitary solutions of this
form to the braid relations based on this list \cite{D}. The
particular importance of one of the solutions,
$\frac{1}{\sqrt{2}}\left(
     \begin{smallmatrix}
       1 & 0 & 0 & 1 \\
       0 & 1 & -1 & 0 \\
       0 & 1 & 1 & 0 \\
       -1 & 0 & 0 & 1 \\
     \end{smallmatrix}
   \right)$, was pointed out in the work of Kauffman and
Lomonaco \cite{KL}, and the connection of R with quantum computing
was explored there.  Additionally, in \cite{FRW} it was shown that
the image of the braid group under the representation afforded by
that particular solution is a finite group.  Furthermore, along with
E. Rowell and Z. Wang, we showed there is an exact sequence:
$$1 \rightarrow E^{-1}_ {n-1} \rightarrow Im(B_n) \rightarrow S_n
\rightarrow 1$$
 $\forall n \geq 2$, where the group $E_m^{-1}$ is
an extraspecial 2 group. Obtaining representations of the braid
group from representations of extraspecial 2-groups also has been
studied in \cite{GJ}.  This paper handles the remaining solutions
given by \cite{D}.\\

One method proposed to build quantum computers is based on topology.
In a topological quantum computer, representations of the braid
group can be used to describe the actions of the quantum bits,
a.k.a. qubits.  It has been conjectured that due to the topological
stability of the system, this construction would be invulnerable to
local errors, which is one of the formidable obstacles to physically
realizing a quantum computer. The qubits themselves are encoded in
the lowest energy states of quasi particles, called anyons, at some
fixed positions in a plane and their use in topological models is
described in \cite{FKLW}. Trajectories of the quasi particles in the
3-dimensional space-time form braids.  Braidings of anyons in a
topological quantum computer change the encoded quantum information.
Representations of braid groups have been proposed as the fractional
statistics of anyons \cite{Wilczek}.

\indent For any given model for quantum computing, it is important
to understand whether or not the model is capable of carrying out
any computation up to any given precision, i.e. whether or not the
model is universal.  In topological models, the universality issue
is translated to a question about the closed images of the braid
group representations.  In particular, for the topological models in
\cite{FKLW}, universality is equivalent to the closure of the images
containing $SU(K_n)$. We do not have density here since we claim:
\begin{thm}
For any unitary $4 \times 4$ solution to the Yang Baxter equation,
$R$, there are representations of the braid group $\pi_R:B_n
\rightarrow U(K_n)$.  If all eigenvalues of $R$ are roots of unity,
there is a short exact sequence
$$1 \rightarrow Im(P_n) \rightarrow Im(B_n) \rightarrow S_n
\rightarrow 1$$ Moreover, $Im(B_n)$ is finite.
\end{thm}
However, there are more elaborate adaptive models which make use of
these cases for quantum computing see \cite{FNW1} \cite{FNW2}.\\

\section{Preliminaries}
Although, any invertible solution to the Yang Baxter equation will
yield a representation of the braid group, with an eye towards
possible applications to quantum computing, we restrict ourselves to
Dye's list since we would need the unitarity condition for physical
realizability.  The following propositions indicate that a single
solution to the Yang Baxter equation produces a class of solutions
by conjugation and scalar multiplication. (See \cite{Kassel}).

It is not sufficient to determine which representations from
\cite{H} are unitary since $ARA^{-1}$ may be unitary when $R$ is
not.
\begin{thm}\cite{D}  There are five families of $4 \times 4$ unitary
matrix solutions to the Yang Baxter equation.  Each has the form
$$kARA^{-1}$$
where k is a scalar with norm 1, $A=Q \otimes Q$, and $Q$ is an
invertible matrix such
that $Q=\left(%
\begin{smallmatrix}
  a & b \\
  c & d \\
\end{smallmatrix}%
\right)$
\begin{enumerate}
\item $R_0=Id$.\\
\item $R_1=\frac{1}{\sqrt{2}}\left(%
\begin{smallmatrix}
  1 & 0 & 0 & 1 \\
  0 & 1 & -1 & 0 \\
  0 & 1 & 1 & 0 \\
  -1 & 0 & 0& 1 \\
\end{smallmatrix}%
\right)$\\  The matrix $Q$ has the following restrictions:
$c=\frac{-a\overline{b}}{d}$ and $|a|=|d|.$
\item $R_2=\left(%
\begin{smallmatrix}
  1 & 0 & 0 & 0 \\
  0 & 0 & \alpha & 0 \\
  0 & \beta & 0 & 0 \\
  0 & 0 & 0 & \gamma \\
\end{smallmatrix}%
\right)$ where
$1=\alpha\overline{\alpha}=\beta\overline{\beta}=\gamma\overline{\gamma}$.\\
The variables in the matrix $Q$ also has the following restriction:
$c=\frac{-a\overline{b}}{d}.$
\item $R_3=\left(%
\begin{smallmatrix}
  0 & 0 & 0 & \alpha \\
  0 & 1 & 0 & 0 \\
  0 & 0 & 1 & 0 \\
  \beta & 0 & 0 & 0\\
\end{smallmatrix}%
\right)$ where $|\alpha\beta|=1$,
$\alpha\overline{\alpha}=\frac{(d\overline{d})^2}{(a\overline{a})^2}$
and
$\beta\overline{\beta}=\frac{(a\overline{a})^2}{(d\overline{d})^2}$.\\
The matrix $Q$ has the restriction: $c=\frac{-a\overline{b}}{d}.$
\item $R_3'=\left(%
\begin{smallmatrix}
  0 & 0 & 0 & \alpha \\
  0 & 1 & 0 & 0 \\
  0 & 0 & 1 & 0 \\
  \beta & 0 & 0 & 0\\
\end{smallmatrix}%
\right)$ where $|\alpha\beta|=1$,
$\alpha=\frac{(b\overline{b}+d\overline{d})(\overline{a}b +
\overline{c}d)}{(a\overline{a}+c\overline{c})(a\overline{b}+c\overline{d})}$
and $\beta=\frac{(a\overline{a}+c\overline{c})(a\overline{b} +
c\overline{d})}{(b\overline{b}+d\overline{d})(\overline{a}b+\overline{c}d)}$.\\
The matrix $Q$ has the restriction: $c \neq
\frac{-a\overline{b}}{d}$.\\
\end{enumerate}
\end{thm}

As the identity case is clearly not interesting and since the
restrictions are of little consequence when looking at images of the
braid groups, for the remainder of this paper we will content
ourselves with analyzing the cases $R_2$ and $R_3$ since
$R_1$ was covered in \cite{FRW}.   We mainly concern ourselves with $$R_2=\left(%
\begin{smallmatrix}
  1 & 0 & 0 & 0 \\
  0 & 0 & \alpha & 0 \\
  0 & \beta & 0 & 0 \\
  0 & 0 & 0 & \gamma \\
\end{smallmatrix}%
\right)$$ since, if we conjugate $R_2$ we get $R_3$ so as abstract
groups, the closed images will be the same.\\
There are common
notations used in
both and we record some here:

Our convention for tensor products of matrices is to use ``left into
right," that is, if $$X=\begin{pmatrix}w & x
\\ y & z\end{pmatrix}~~ and ~~A=\begin{pmatrix}a & b \\ c &
d\end{pmatrix}~~then~~X\otimes A=\begin{pmatrix} aX & bX \\ cX &
dX\end{pmatrix}.$$

We let $$\pi_n(\sigma_i)=I_2^{\otimes i-1}\otimes R\otimes
I_2^{\otimes n-i-1}$$ be the representations of the braid group,
$B_n$ arising from $R$.

\section{The Image of the Braid Group}

\subsection{Restriction to $P_n$}

Notice that $$R=
\left(%
\begin{smallmatrix}
  1 & 0 & 0 & 0 \\
  0 & 0 & \alpha & 0 \\
  0 & \beta & 0 & 0 \\
  0 & 0 & 0 & \gamma \\
\end{smallmatrix}%
\right) =
\left(%
\begin{smallmatrix}
  1 & 0 & 0 & 0 \\
  0 & 0 & 1 & 0 \\
  0 & 1 & 0 & 0 \\
  0 & 0 & 0 & 1 \\
\end{smallmatrix}%
\right)
\left(%
\begin{smallmatrix}
  1 & 0 & 0 & 0 \\
  0 & \beta & 0 & 0 \\
  0 & 0 & \alpha & 0 \\
  0 & 0 & 0 & \gamma \\
\end{smallmatrix}%
\right) =P \cdot D$$\\
Where $P$ is a permutation matrix and $D$ is a diagonal matrix.
\begin{prop} \label{abelian}
$Image(P_n)$ is abelian
\begin{proof}
First notice that the image of $\sigma_i^2$ is a diagonal matrix for
each $i$.  Now consider
\begin{eqnarray*}
\lefteqn{\pi_n(\sigma_{i+1}) \pi_n(\sigma_i^2) \pi_n(\sigma_{i+1}^{-1})}\\
&=& (I^{\otimes i} \otimes R \otimes I^{\otimes n-i-2})(I^{\otimes
i-1} \otimes R^2 \otimes I^{\otimes n-i-1})(I^{\otimes i} \otimes R
\otimes I^{\otimes n-i-2})^{-1}\\
 &=& I^{\otimes i-1}[(I \otimes R)(R^2 \otimes I)(I \otimes R^{-1})]
 \otimes I^{\otimes n-i-2}\\
 &=& I^{\otimes i-1} \otimes [(I \otimes PD)(R^2 \otimes I)(I
 \otimes(PD)^{-1})] \otimes I^{\otimes n-i-2}\\
 &=& I^{\otimes i-1} \otimes [(I \otimes PD)(R^2 \otimes I)(I
 \otimes D^{-1}P)] \otimes I^{\otimes n-i-2}\\
 &=& I^{\otimes i-1} \otimes [(I \otimes P)(I \otimes D)(R^2 \otimes I)(I
 \otimes D^{-1})(I \otimes P)] \otimes I^{\otimes n-i-2}\\
 &=& I^{\otimes i-1} \otimes [(I \otimes P)(R^2 \otimes I)(I
 \otimes P)] \otimes I^{\otimes n-i-2}\\
 &=& D'
\end{eqnarray*}
where $D'$ is a diagonal matrix.  We only used that $R^2 \otimes I$,
$I \otimes D$ and $D \otimes I$ are a diagonal matrices and
therefore commute. Since the pure braid group, $P_n$, is generated
by all conjugates of $\sigma_i^2$ the above shows that the
Image$(P_n)$ is a subset of the diagonal matrices and thus abelian.
\end{proof}
\end{prop}

\begin{rmk}
Since $R$ has arbitrary variables, it might be possible for $R$
itself to generate an infinite group.  So for the remainder of this
paper we restrict ourselves to the case where the eigenvalues of $R$
are roots of unity and thus $\exists k$ such that $p_i^k=I$, where
$p_i=\pi_n((\sigma_i)^2)=I_2^{\otimes (i-1)} \otimes R^2 \otimes
I_2^{\otimes (n-i-1)}$. This in combination with the preceding
proposition allow us to conclude the following:
\end{rmk}

\begin{prop}
$Image(P_n)$ is a finite abelian group.
\end{prop}

\begin{theorem}\label{quotient2}
We have an exact sequence:
$$1\rightarrow Im(P_n) \rightarrow Im(B_n) \rightarrow S_n\rightarrow 1$$ for all $n\geq 3$.
In other words, $G_n$ is an extension of $H_n$ by $S_n$.
\end{theorem}
\begin{proof}
Notice $\pi_n(B_n)/\pi_n(P_n)$ is a homomorphic image of $S_n$ as
$\pi_n$ induces a surjective homomorphism $\Hat{\pi}_n:
B_n/P_n\rightarrow \pi_n(B_n)/\pi_n(P_n)$ and $B_n/P_n\cong S_n$. We
would like to know if $\Hat{\pi}_n$ is an isomorphism in this case
as well.  Therefore, we must determine if $\Ker(\Hat{\pi}_n)$ is
trivial. We note for $n\geq 4$ it is sufficient to check that the
element $(12)(34)$ is not in the kernel, while for $n=3$ we should
check that $(123)$ is not in the kernel.  We simply observe the
corresponding elements $\sigma_1\sigma_3$ and $\sigma_2\sigma_1$ are
not diagonal and note that in Proposition \ref{abelian} we showed
the image of the pure braid group is a subset of diagonal matrices.
\end{proof}
In particular, we have shown that the image of the braid group is
finite as long as the eigenvalues are roots of unity.

\subsection{Invariants} We now put our representation through the
machinery given by Turaev \cite{Turaev} and uncover link invariants
if
and only if $\gamma =\pm 1$ where: $$R=\left(%
\begin{smallmatrix}
  1 & 0 & 0 & 0 \\
  0 & 0 & \alpha & 0 \\
  0 & \beta & 0 & 0 \\
  0 & 0 & 0 & \gamma \\
\end{smallmatrix}%
\right).$$  If $\gamma=1$, then $\mu=xy\cdot Id$ and the invariant
is:
$$T_{R,x}=x^{n-e(\sigma)}\textrm{Trace}(\pi_n(\sigma)).$$

\noindent If $\gamma=-1$, then $\mu=xy\left(%
\begin{smallmatrix}
  1 & 0 \\
  0 & -1 \\
\end{smallmatrix}%
\right)$ and the invariant is:
$$T_{R,x}=x^{n-e(\sigma)}\textrm{Trace}(\mu \otimes \mu \circ \pi_n(\sigma)).$$

We would like to know if these invariants correspond to some known
invariant. We proceed by breaking into cases with the number of
eigenvalues of $R$.  Note that if $R$ has only 2 eigenvalues they
must be $\pm 1$ and we get $\pi_n(B_n) \cong S_n$.
 If we let $\gamma=1$ then the other eigenvalues of $R$ are
 $\pm \sqrt{\alpha\beta}$ so we will have exactly 3 eigenvalues if $\alpha\beta \neq 1$.  If
 we let $\gamma=-1$ and assume that $\alpha\beta \neq 1$ we will
 have 4 distinct eigenvalues.

\section{Different Eigenvalues}

 \subsection{Three eigenvalues}
Let $A_n$ denote the group algebra of the image of the braid group
afforded by: $$R=\left(%
\begin{smallmatrix}
  1 & 0 & 0 & 0 \\
  0 & 0 & \alpha & 0 \\
  0 & \beta & 0 & 0 \\
  0 & 0 & 0 & 1 \\
\end{smallmatrix}%
\right)$$  The main result of what follows will hold for $R_3$ and
the arguments are analogous.  Recall that the BMW algebras,
$C_n(r,q)$, are complex group algebras with generators $g_i$
satisfying the braid relations and
\begin{enumerate}
\item[(R1)] $e_ig_i=r^{-1}e_i$
        \item[(R2)] $e_ig_{i-1}^{\pm 1}e_i=r^{\pm 1} e_i$
    \end{enumerate}
where $e_i=1-\frac{(g_i-g_i^{-1})}{(q-q^{-1})}$.\\
\indent
These are quotients of the braid group and were originally defined
by Birman and Wenzl \cite{BW} and independently by Murakami
\cite{M}. However, we are following the definition given in
\cite{Wen} where Wenzl shows if $r = \pm q^k$, $C_n(r,q)$ is not semisimple, but there is a
trace on this algebra (corresponding to the Kauffman
polynomial) and moding out by the annihilator of this trace recovers
semisimplicity \cite{Wen}. \\

 We wish to make a
connection between these algebras and $A_n$ so we take the following
steps: \begin{enumerate} \item We rescale our matrix $R$ and by
abuse of notation let $R=QR$ \item We next suppose that
$\alpha,~\beta=Q^{-2}$.
\end{enumerate}
  Then $$R=\left(%
\begin{smallmatrix}
  Q & 0 & 0 & 0 \\
  0 & 0 & Q^{-1} & 0 \\
  0 & Q^{-1} & 0 & 0 \\
  0 & 0 & 0 & Q \\
\end{smallmatrix}%
\right)$$  Now we define a map $\phi:C_n(q,q) \rightarrow A_n$.  To
avoid confusion, we denote elements in $C_n(q,q)$ as lower case
letters and their image in capital letters: $\phi(g_i)=G_i$ where
$G_i=\pi_n(\sigma_i)$ and $q=Q$.  Notice that $G_1=R$ and and we have a nice form for $$E_i=I^{\otimes (i-1)} \otimes \left(%
\begin{smallmatrix}
  0 & 0 & 0 & 0 \\
  0 & 1 & 1 & 0 \\
  0 & 1 & 1 & 0 \\
  0 & 0 & 0 & 0 \\
\end{smallmatrix}%
\right) \otimes I^{\otimes (n-i-1)}$$

\begin{lemma}The $G_i$ 's satisfy
$(R1)$ and $(R2)$ so $\phi$ is a homomorphism.
\begin{proof}
For (R1) notice that
\begin{align*}
E_iG_i\\
& =\bigl(I^{\otimes (i-1)} \otimes \left(%
\begin{smallmatrix}
 0 & 0 & 0 & 0 \\
 0 & 1 & 1 & 0 \\
 0 & 1 & 1 & 0 \\
 0 & 0 & 0 & 0 \\
\end{smallmatrix}%
\right) \otimes I^{\otimes (n-i-1)}\bigr)
\bigl(I^{\otimes (i-1)} \otimes \left(%
\begin{smallmatrix}
 Q&0&0&0 \\
0 & 0 & Q^{-1} & 0 \\
0 & Q^{-1} & 0 & 0 \\
0 & 0 & 0 & Q \\
\end{smallmatrix}%
\right) \otimes I^{\otimes (n-i-1)}\bigl)\\
&=I^{\otimes (i-1)} \otimes \left(%
\begin{smallmatrix}
  Q & 0 & 0 & 0 \\
  0 & 0 & Q^{-1} & 0 \\
  0 & Q^{-1} & 0 & 0 \\
  0 & 0 & 0 & Q \\
\end{smallmatrix}%
\right) \left(%
\begin{smallmatrix}
  0 & 0 & 0 & 0 \\
  0 & 1 & 1 & 0 \\
  0 & 1 & 1 & 0 \\
  0 & 0 & 0 & 0 \\
\end{smallmatrix}%
\right) \otimes I^{\otimes (n-i-1)}\\
&=I^{\otimes (i-1)} \otimes  \left(%
\begin{smallmatrix}
  0 & 0 & 0 & 0 \\
  0 & Q^{-1} & Q^{-1} & 0 \\
  0 & Q^{-1} & Q^{-1} & 0 \\
  0 & 0 & 0 & 0 \\
\end{smallmatrix}%
\right) \otimes I^{\otimes (n-i-1)}\\
&=Q^{-1}(I^{\otimes (i-1)} \otimes \left(%
\begin{smallmatrix}
  0 & 0 & 0 & 0 \\
  0 & 1 & 1 & 0 \\
  0 & 1 & 1 & 0 \\
  0 & 0 & 0 & 0 \\
\end{smallmatrix}%
\right) \otimes I^{\otimes (n-i-1)})\\
&=Q^{-1}E_i
\end{align*}

 For (R2) notice that

\begin{align*}
E_iG_{i-1}E_i\\
&= \left(I^{\otimes (i-1)} \otimes \left(%
\begin{smallmatrix}
  0 & 0 & 0 & 0 \\
  0 & 1 & 1 & 0 \\
  0 & 1 & 1 & 0 \\
  0 & 0 & 0 & 0 \\
\end{smallmatrix}%
\right) \otimes I^{\otimes (n-i-1)}\right)(I^{\otimes ((i-1)-1)}
\otimes R \otimes I^{(n-i)})\\
& \left(I^{\otimes (i-1)} \otimes \left(%
\begin{smallmatrix}
  0 & 0 & 0 & 0 \\
  0 & 1 & 1 & 0 \\
  0 & 1 & 1 & 0 \\
  0 & 0 & 0 & 0 \\
\end{smallmatrix}%
\right) \otimes I^{\otimes (n-i-1)}\right)\\
&=\left(I^{\otimes (i-2)} \otimes \left(I \otimes \left(%
\begin{smallmatrix}
  0 & 0 & 0 & 0 \\
  0 & 1 & 1 & 0 \\
  0 & 1 & 1 & 0 \\
  0 & 0 & 0 & 0 \\
\end{smallmatrix}%
\right)\right) \otimes I^{\otimes (n-i-1)}\right)((I^{\otimes (i-2)}
\otimes (R \otimes I)) \otimes I^{(n-i-1)})\\
& ~\left(I^{\otimes (i-2)} \otimes \left(I \otimes \left(%
\begin{smallmatrix}
  0 & 0 & 0 & 0 \\
  0 & 1 & 1 & 0 \\
  0 & 1 & 1 & 0 \\
  0 & 0 & 0 & 0 \\
\end{smallmatrix}%
\right) \otimes I^{\otimes (n-i-1)}\right)\right)\\
&=I^{\otimes (i-2)}\otimes \left(I \otimes \left(%
\begin{smallmatrix}
  0 & 0 & 0 & 0 \\
  0 & 1 & 1 & 0 \\
  0 & 1 & 1 & 0 \\
  0 & 0 & 0 & 0 \\
\end{smallmatrix}%
\right)\right)(R \otimes I)\left(I \otimes \left(%
\begin{smallmatrix}
  0 & 0 & 0 & 0 \\
  0 & 1 & 1 & 0 \\
  0 & 1 & 1 & 0 \\
  0 & 0 & 0 & 0 \\
\end{smallmatrix}%
\right)\right) \otimes I^{\otimes (n-i-1)}\\
&=I^{\otimes (i-2)}\otimes \left(%
\begin{smallmatrix}
  0 & 0 & 0 & 0 & 0 & 0 & 0 & 0 \\
  0 & 0 & 0 & 0 & 0 & 0 & 0 & 0 \\
  0 & 0 & Q & 0 & Q & 0 & 0 & 0 \\
  0 & 0 & 0 & Q & 0 & Q & 0 & 0 \\
  0 & 0 & Q & 0 & Q & 0 & 0 & 0 \\
  0 & 0 & 0 & Q & 0 & Q & 0 & 0 \\
  0 & 0 & 0 & 0 & 0 & 0 & 0 & 0 \\
  0 & 0 & 0 & 0 & 0 & 0 & 0 & 0 \\
\end{smallmatrix}%
\right)\otimes I^{\otimes(n-i-1)}\\
&=I^{\otimes(i-2)}\otimes \left(I \otimes \left(%
\begin{smallmatrix}
  0 & 0 & 0 & 0 \\
  0 & Q & Q & 0 \\
  0 & Q & Q & 0 \\
  0 & 0 & 0 & 0 \\
\end{smallmatrix}%
\right)\right) \otimes I^{\otimes (n-i-1)}\\
&=Q I^{\otimes(i-1)} \otimes \left(%
\begin{smallmatrix}
  0 & 0 & 0 & 0 \\
  0 & 1 & 1 & 0 \\
  0 & 1 & 1 & 0 \\
  0 & 0 & 0 & 0 \\
\end{smallmatrix}%
\right)\otimes I^{\otimes (n-i-1)}\\
& = Q \cdot E_i
\end{align*}
\end{proof}
\end{lemma}

So we have:
\[
\begin{diagram}
\node{C_n(q,q)}
      \arrow[2]{e,t}{\phi}
   \node[2]{A_n}
\\
\node[2]{B_n} \arrow{ne,t}{\pi_n} \arrow{nw,t}{}
\end{diagram}
\]
Where the map from $B_n \rightarrow C_n(q,q)$ is given my $\sigma_i
\mapsto g_i$.\\

\indent We also note here that by taking the standard trace of the
matrices $E_i$ and $G_i$ we have that
$\textrm{Trace}(E_i)=\frac{1}{2} \cdot 2^{n}$ and
$\textrm{Trace}(G_i)=Q \cdot \frac{1}{2} \cdot2^{n}.$

\begin{prop}\cite{Wen} Let $x=1+ \frac{r-r^{-1}}{q-q^{-1}}$.  There
exists a functional trace, $tr$, on $C_\infty(r,q)$ \emph{uniquely}
defined inductively by:
\begin{enumerate}
\item $tr(1)=1$
\item $tr(ab)=tr(ba)$
\item $tr(e_i)=\frac{1}{x}$
\item $tr(g_i^{\pm 1})=\frac{r^{\pm 1}}{x}$
\item $tr(a \chi b)=tr(\chi)tr(ab)$ for $a,b \in C_{n-1}(r,q)$,
$\chi=g_{n-1}~or~e_{n-1}$
\end{enumerate}
\end{prop}

We are in the case $r=q$ and so $x=2$.  Note that the standard Trace
on matrices does not satisfy these conditions, but if we let
$Tr_n=\frac{1}{2^{n}}\textrm{Trace}_n$ where Trace is the standard
trace, we have all but the last condition.  That is, we wish to show
$Tr(A \chi B)=Tr(\chi)Tr(AB)$ for $\chi \in \{E_n,~G_n\}$ and $A,~B
\in A_n$.  Note that $A_n \hookrightarrow A_{n+1}$ by $X \mapsto X
\otimes I_2$, so since $A,~B~ \in A_n$ and $\chi \in
\{E_n,~G_n\}\subset A_{n+1}$, we actually consider
$A\otimes I_2$ and $B\otimes I_2$ so that $A\chi B$ makes sense.  Let $\chi=G_n=I^{\otimes n-1} \otimes R$.
Then we have:\\

\begin{align*}
Tr_{n+1}(A \chi B) \\
&= Tr_{n+1}((A \otimes I) \chi (B \otimes I)) \\
& = Tr_{n+1}((A \otimes I)(I^{\otimes n-1} \otimes R)(B \otimes I))\\
& = Tr_{n+1}(I^{\otimes n-1} \otimes R)(B \otimes I)(A \otimes I))\\
& = Tr_{n+1}\left((I^{\otimes n - 1} \otimes R)
\left(%
\begin{smallmatrix}
  BA & 0_{2^n} \\
  0_{2^n} & BA \\
\end{smallmatrix}%
\right)\right)\\
&= Tr_{n+1}\left(\left(%
\begin{smallmatrix}
  \left(%
\begin{smallmatrix}
  Q \cdot I^{\otimes n-1} & 0_{2^n -1} \\
  0_{2^n -1} & 0_{2^n -1} \\
\end{smallmatrix}%
\right) & \left(%
\begin{smallmatrix}
  0_{2^n -1} & 0_{2^n -1} \\
  Q \cdot \alpha \cdot I^{\otimes n -1} & 0_{2^n -1} \\
\end{smallmatrix}%
\right) \\
  \left(%
\begin{smallmatrix}
  0_{2^n -1} & Q \cdot \beta \cdot I^{\otimes n -1} \\
  0_{2^n -1} & 0_{2^n -1} \\
\end{smallmatrix}%
\right) & \left(%
\begin{smallmatrix}
  0_{2^n -1} & 0_{2^n -1} \\
  0_{2^n -1} & Q \cdot I^{\otimes n-1} \\
\end{smallmatrix}%
\right) \\
\end{smallmatrix}%
\right) \left(%
\begin{smallmatrix}
  BA & 0_{2^n} \\
  0_{2^n} & BA \\
\end{smallmatrix}%
\right)\right) \\
& = Tr_{n+1}\left(\left(%
\begin{smallmatrix}
  \left(%
\begin{smallmatrix}
  Q \cdot I^{\otimes n-1} & 0_{2^n -1} \\
  0_{2^n -1} & 0_{2^n -1} \\
\end{smallmatrix}%
\right)BA & \left(%
\begin{smallmatrix}
  0_{2^n -1} & 0_{2^n -1} \\
  Q \cdot \alpha \cdot I^{\otimes n -1} & 0_{2^n -1} \\
\end{smallmatrix}%
\right)BA \\
  \left(%
\begin{smallmatrix}
  0_{2^n -1} & Q \cdot \beta \cdot I^{\otimes n -1} \\
  0_{2^n -1} & 0_{2^n -1} \\
\end{smallmatrix}%
\right)BA & \left(%
\begin{smallmatrix}
  0_{2^n -1} & 0_{2^n -1} \\
  0_{2^n -1} & Q \cdot I^{\otimes n-1} \\
\end{smallmatrix}%
\right)BA \\
\end{smallmatrix}%
\right)\right)\\
&=Tr_n\left(\left(%
\begin{smallmatrix}
  Q \cdot Id^{\otimes n-1} & 0_{2^n -1} \\
  0_{2^n -1} & 0_{2^n -1} \\
\end{smallmatrix}%
\right)BA\right) + Tr_n \left(\left(%
\begin{smallmatrix}
   0_{2^n -1} & 0_{2^n -1} \\
  0_{2^n -1} & Q \cdot I^{\otimes n-1} \\
\end{smallmatrix}%
\right)BA\right)\\
&=Tr_n\left(\left(%
\begin{smallmatrix}
  Q \cdot Id^{\otimes n-1} & 0_{2^n -1} \\
  0_{2^n -1} & 0_{2^n -1} \\
\end{smallmatrix}%
\right)BA + \left(%
\begin{smallmatrix}
   0_{2^n -1} & 0_{2^n -1} \\
  0_{2^n -1} & Q \cdot I^{\otimes n-1} \\
\end{smallmatrix}%
\right)BA\right)\\
&=Tr_n\left(\left(\left(%
\begin{smallmatrix}
  Q \cdot I^{\otimes n-1} & 0_{2^n -1} \\
  0_{2^n -1} & 0_{2^n -1} \\
\end{smallmatrix}%
\right) + \left(%
\begin{smallmatrix}
   0_{2^n -1} & 0_{2^n -1} \\
  0_{2^n -1} & Q \cdot I^{\otimes n-1} \\
\end{smallmatrix}%
\right)\right)BA\right)\\
&= Tr_n\left(\left(%
\begin{smallmatrix}
  Q \cdot I^{\otimes n-1} & 0_{2^n -1} \\
  0_{2^n -1} & Q \cdot I^{\otimes n -1} \\
\end{smallmatrix}%
\right)BA\right)\\
&=Tr_n (Q \cdot I^{\otimes n} \cdot BA) \\
&= Q \cdot Tr_n (BA)\\
&= \frac{1}{2} Q \cdot Tr_{n+1} (BA)\\
&= \frac{1}{2} \cdot Tr_{n+1} (AB)\\
&= Tr_{n+1}(G_n) \cdot Tr_{n+1}(AB)
\end{align*}

We have shown the case that $\chi=G_n$, but the proof for $\chi=E_n$
is virtually the same and left to the reader.  Therefore, we have
that the traces correspond.

Since the trace on $C_n(q,q)$ computes the 2 variable Kauffman
polynomial, we have shown in the 3 eigenvalue case if we let
$\alpha, \beta =q^{-2}$ and rescale our representation, the
invariant specified by Turaev is related to a specific
specialization of the 2 variable Kauffman polynomial at $q$.

\begin{prop}\cite{Wen}
$\overline{C_n(q,q)}:=C_n(q,q)/Ann(q,q)$ is semisimple\\ where
$Ann(q,q)$ is the annihilator of the trace.  $Ann(q,q)=\{b \in
C_n(q,q)~|~tr(ab)=0 ~\forall a \in C_n(q,q)\}$
\end{prop}

\begin{lemma}
The induced map $$\overline{\phi}:\overline{C_n(q,q)} \rightarrow
\phi(C_n(q,q))/\phi(Ann(q,q))$$ is injective.
\end{lemma}
\begin{proof}
It suffices to show $ker(\phi)\subset Ann(q,q)$.  But since we have
shown $tr$ corresponds to $Tr$, this is not difficult.  Let $a \in
ker(\phi)$ and $b \in C_n(q,q)$, then $tr(ab)=Tr(ab)=0$, that is, $a
\in Ann(q,q)$
\end{proof}
We note the above proof was inspired by \cite{LR}.\\
An immediate corollary is that $\overline{C_n(q,q)}$ is isomorphic
to a quotient of $A_n$ which we know is finite and thus
$\overline{C_n(q,q)}$ is finite.

\begin{remark} The question of identifying the closed image of the braid
groups in unitary representations associated with the 2 variable
Kauffman polynomial at specific roots of unity was begun in
\cite{LRW} and continued in \cite{LR}.  The focus of \cite{LRW} is
on the cases where the closed images are infinite; in particular,
they contain the special unitary group. There are open questions
about the cases where the images are finite.  The BMW algebras are
indexed by $q$ and $r$ and when the image of the braid group
representations are known to be finite $q$ and $r$ are related in
some way.  The case $r=q$ was covered above which only leaves the
case $r=q^{\pm \ell/2}$.
\end{remark}

 \subsection{Four eigenvalues}
In \cite{FRW} it was shown that after normalization certain
solutions of the matrices from \cite{D} are related to the Jones
polynomial which has a skein relation where $\sigma_i$ is the sum of
2 terms. This paper has shown similar results for the Kauffman
polynomial which has a skein relation where $\sigma_i$ is the sum of
3 terms. We now note that the $G_2$ invariant defined by Kuperberg
\cite{Kuperberg} has a skein relation in which $\sigma_i$ is the sum
of 4 terms. When $\gamma = -1$ by \cite{Turaev} we get an invariant
and with the aid of computers we can show in some cases our matrices
correspond to particular specializations the $G_2$ invariant. We
list them now for completeness:
\begin{thm}
Let $x=\sqrt{\alpha\beta}$ and let $t$ be the variable in
Kuperberg's invariant.  Then for the following specializations of
$x$ and $t$ our matrices satisfy Kuperberg's skein relation:
\begin{enumerate}
\item $t=1$, $x$ arbitrary
\item $t=i,~x=-i$ and $t=-i,~x=i$ and $t=\pm i,~x=\pm i$
\item $t=\frac{1}{2}(1 \pm \sqrt{3})$, $x=-t$ and $x=t-1$
\end{enumerate}
\end{thm}
These were found by using the skein relation as follows:\\
First, we noted that the minimal polynomial of $R$ gave us:
\[(R^2 - \alpha\beta)(R^2-1)=0\]
\[\Rightarrow R^4-(\alpha\beta +1)R^2=\alpha\beta=0\]
\[\Rightarrow R^4=(\alpha\beta+1)R^2 -\alpha\beta\]
Then we identified $R$ with a crossing and using Kuperberg's skein
relation resolved each crossing.  Finally, we used a computer to
equate like coefficients and found the preceding solutions.  These
were found using Maple and checked by the author by hand, but
further calculations have proven fruitless.

\end{document}